\documentclass[12pt]{amsart}

\usepackage{amsfonts,amsmath,latexsym,amssymb,verbatim,amsbsy}
\usepackage{amsthm, pdfsync}

\numberwithin{equation}{section}

\theoremstyle{plain}
\newtheorem{THEOREM}{Theorem}[section]

\theoremstyle{definition}

\theoremstyle{remark}

\newtheorem{remark}[THEOREM]{Remark}

\def \a {\alpha}

\def \d {\delta}
\def \g {\gamma}
\def \e {\varepsilon}

\def \k {\kappa}

\def \n {\nabla}

\def \p {\partial}
\def \ra {\rightarrow}

\newcommand{\der}[2]{(#1 \cdot \triangledown) #2}

\newcommand{\ddt}{\frac{d}{dt}}

\DeclareMathOperator{\supp}{supp} %

\begin{document}

\title{A continuous model for turbulent energy cascade}

\author{A. Cheskidov}
\address[A. Cheskidov, R. Shvydkoy]
{Department of Mathematics, Stat. and Comp. Sci.\\
       University of Illinois\\
       Chicago, IL 60607}
\email{acheskid@math.uic.edu}

\author{R. Shvydkoy}

\email{shvydkoy@math.uic.edu}

\author{S. Friedlander}
\address[S. Friedlander]{Department of Mathematics\\
University of Southern California\\
3620 S. Vermont Ave.\\
Los Angeles, CA 90089}
\email{susanfri@usc.edu}

\thanks{The work of A. Cheskidov is partially supported by NSF grant DMS--0807827, S. Friedlander is partially supported by the NSF grant DMS--0849397, R. Shvydkoy acknowledges the support of NSF grant DMS-- 0907812.}

\begin{abstract}
In this paper we introduce a new PDE model in frequency space for the inertial energy cascade that reproduces the classical scaling laws of Kolmogorov's theory of turbulence. 
Our point of view is based upon studying the energy flux through a continuous range of scales rather than the discrete set of dyadic scales. The resulting model is a variant of Burgers equation on the half line with a boundary condition which represents a constant energy input at integral scales. The viscous dissipation is modeled via a damping term. We show existence of a unique stationary solution, both in the viscous and inviscid cases, which replicates the classical dissipation anomaly in the limit of vanishing viscosity. 

A survey of recent developments in the deterministic approach to the laws of turbulence, and in particular, to Onsager's conjecture is given.
\end{abstract}



\maketitle

\section{Motivation for the Model}

\subsection{Onsager and Kolmogorov}

The Euler equations for the motion of an incompressible,
inviscid fluid are
\begin{equation}\label{e:ee}
\frac{\p u}{\p t} + \der{u}{u} = - \n p
\end{equation}
\begin{equation}\label{e:div}
\n \cdot u = 0,
\end{equation}
where $u(x, t)$ is a divergence free velocity vector and $p(x, t)$
is the internal pressure.  We consider the system in 3 spatial
dimensions and we assume that the domain is either periodic or the
entire space $\Bbb{R}^3$. To obtain the energy equation we multiply
\eqref{e:ee} by $u$ and integrate, using \eqref{e:div} to give
\begin{equation}\label{e:balance}
\frac{1}{2} \frac{d}{dt} \int u^2 dx = -\int (u \cdot
\triangledown) u \cdot u dx.
\end{equation}
We define the total energy flux $\Pi$ by
\begin{equation}\label{e:flux}
\Pi = \int (u \cdot \triangledown) u \cdot u dx. 
\end{equation}
For {\it smooth} solutions we can integrate by parts and use \eqref{e:div}
to conclude that $\Pi = 0$ and hence energy conservation holds, i.e.
\begin{equation}\label{e:ene}
\int |u(x, t)|^2 dx = \int |u(x, 0)|^2 dx \qquad \mbox{ for } t \ge 0.
\end{equation}
However, in the context of turbulent flows in the limit of vanishing
viscosity, it is appropriate to consider the Euler equations in the
sense of distributions and impose only minimal assumptions on the
regularity of the velocity field $u$.  In the absence of sufficient
smoothness we cannot integrate by parts in \eqref{e:flux} or even make sense of \eqref{e:flux} and ensure that $\Pi =
0$.  Conservation of energy might then be violated.  Hence it is of
interest to ask what are the minimal regularity assumptions on the
velocity that ensures that \eqref{e:ene} holds.

Observing that the integrand in \eqref{e:flux} is cubic in $u$ and contains
one spatial derivative suggests that if $u$ has Holder continuity $h
> 1/3$, integration by parts is justified and $\Pi = 0$.  In fact
this was the conjecture made many years ago by Onsager in his
seminal paper on statistical fluid dynamics \cite{onsager}. More
precisely, he conjectured that (a) every weak solution to the Euler
equation with smoothness $h > 1/3$ conserves energy and (b) there
exists a weak solution with $h \le 1/3$ which does not conserve
energy. Such putative energy dissipation due to the irregularity of
the flow is called anomalous or turbulent dissipation.  A detailed
historical account of Onsager's theory is given by Eyink and
Sreenivasan \cite{es}.

All physical fluids are viscous, if only very weakly so.  Turbulent
fluids are believed to be described by the Navier-Stokes equations
\begin{align}
\frac{\partial u}{\partial t} + (u \cdot \triangledown) u & = -
\triangledown P + \nu \triangle u + f  \label{e:nse}\\
\triangledown \cdot u &= 0, \label{e:nsediv}
\end{align}
where $\nu$, which could be very small, is the coefficient of
viscosity and $f$ is an external force which supplies energy into
the system. The ``classical" Kolmogorov theory of turbulence
predicts that energy dissipative solutions to the Euler equation may
arise in the limit of vanishing viscosity for ``generic" viscous
flows that are governed by \eqref{e:nse} - \eqref{e:nsediv}. In homogeneous, isotropic
turbulence the mean kinetic energy per unit mass is defined by
$E  = \frac{1}{2} \langle|u|^2\rangle$ while the energy
density spectrum is defined by $E(\k) = \frac{1}{2} \frac{d}{d\k}
\langle|u_{<\k}|^2\rangle$. Here $u_{<\k}$ denotes the filtered
velocity field containing all frequencies below a wave number $\k$.
Hence $E = \int^\infty_0 E(\k)d\k$. The mean energy
dissipation rate per unit mass is defined by
\begin{equation}\label{e:dissrate}
\varepsilon^\nu = \langle \nu|\triangledown u^\nu|^2\rangle,
\end{equation}
where $u^\nu$ is a solution to \eqref{e:nse} - \eqref{e:nsediv}. Kolmogorov
\cite{K41} predicted that the energy cascade mechanism in fully
developed 3-dimensional turbulence produces a striking phenomenon,
namely the persistence of non-vanishing energy dissipation in the
limit of vanishing viscosity, i.e.
\begin{equation}\label{e:Kol}
\lim_{\nu\to 0} \varepsilon^\nu \to \varepsilon > 0
\end{equation}
where $\varepsilon$ is the anomalous dissipation rate for the
inviscid Euler system. The positivity of the limit in \eqref{e:Kol} is
referred to as the dissipation anomaly.  Let us now
assume that $f = f_{<\k_{f}}$ (i.e. $f$ has finite Fourier support)
and that solutions to \eqref{e:nse} tend to a statistically stationary state
with uniformly bounded mean energy.  We multiply \eqref{e:nse} by
$u^\nu_{<\k}$ and obtain
\begin{equation}\label{e:fluxnse}
\Pi_\nu (\k) = - \nu \langle|\triangledown u^\nu_{<\k}|^2\rangle + \langle
f \cdot u^\nu_{<\k}\rangle.
\end{equation}
If $\k > \k_f$ we have $\langle f \cdot u^\nu_{<\k}\rangle = \langle f
\cdot u^\nu\rangle = \varepsilon^\nu > 0$.  On the other hand, by
Bernstein's inequality, $\nu\langle|\triangledown
u^\nu_{<\k}|^2\rangle \le \nu \k^2 \langle|u^\nu|^2\rangle$. Since the
energy is uniformly bounded by assumption, we obtain from \eqref{e:fluxnse}
that
\begin{equation}\label{e:fluxvan}
\lim_{\nu \to 0} \Pi_\nu(\k) = \lim_{\nu \to 0} \varepsilon^\nu =
\varepsilon.
\end{equation}
Thus in the limit of vanishing viscosity the average solution of the
forced Euler equation inherits the anomalous dissipation rate
$\varepsilon$.

As Frisch \cite{frisch} describes, a self-similarity hypothesis on
the velocity increments in small (spatial) scales implies that the
energy spectrum as a function of wave number $\k$ has the power law
\begin{equation}\label{e:powerlaw}
E(\k) \sim \varepsilon^{2/3} \k^{-5/3} 
\end{equation}
in the ``inertial" range $\k \in[\k_f, \k_d]$. Here $\k_d$ is the
Kolmogorov dissipation wave number given by $\k_d =
(\varepsilon/\nu^3)^{1/4}$ and $\k_f = \max \{|\k|: \k \in
\supp \hat{f}\}$.  This power law is known as the $K41$
turbulence model. Although the 5/3 power law is consistent with much
physical data, there are also experiments which indicate turbulent
regimes with alternative power laws.  In fact, Kolmogorov's 1941
theory requires that the local velocity fluctuations are uniformly
distributed over space.  However, in reality dynamical stretching of
the vortex filaments in 3-dimensional flows leaves some regions of
the fluid domain with moderate turbulent activity and other regions
with intense activity.  This so called spatial intermittency should
reasonably be accounted for  in the description of the scaling laws.
The expressions for $E(\k)$ and $\k_d$ that incorporate the dimension
$D$ of the effective dissipation region are
\begin{equation}\label{e:specinter}
E(\k) \sim \varepsilon^{2/3} \k^{-(8-D)/3}
\end{equation}
and
\begin{equation}\label{e:dissinter}
\k_d \sim (\varepsilon/\nu^3)^{1/1+D} 
\end{equation}
where $D \in [0, 3]$. Thus the classical $K41$ model corresponds to
$D=3$, i.e. uniform distribution over 3 dimensional space, while
$D=0$ corresponds to a fully intermittent model where energy
cascades through scales and dissipates only on points.

\subsection{Onsager's Conjecture and Besov Spaces}

In the past few years there have been a number of articles that
address part (a) of Onsager's conjecture.  These include articles by
Constantin et al \cite{cet}, Eyink \cite{e}, Duchon and
Robert \cite{dr}. It was shown that appropriate function
spaces to examine the Euler equations in the context of Onsager's
conjecture are Besov spaces.  In such spaces the notion of energy
balance when the velocity is ``a little smoother" than Holder $h >
1/3$ can be made precise. These are the natural spaces to work with
in terms of a description of the energy flux phrased by a
Littlewood-Paley decomposition which provides detailed information
concerning the cascade of energy.  Recently Cheskidov et al
\cite{ccfs} obtained the largest Besov space where conservation
of energy is ensured for the Euler equation. We note that to date
there are no examples of Euler flows that possess some smoothness
and confirm the second part of Onsager's conjecture, although there
are examples of ``very weak" Euler solutions that violate the energy
balance condition \cite{ds,scheffer,shn}.

We recall the definition of a weak solution of the Euler equation. A vector field $u \in C_w([0, T]:
L^2(\Bbb{R}^3))$ is a weak solution of the Euler equations with
initial data $u_0 \in L^2(\Bbb{R}^3)$ if for every compactly
supported test function $\psi \in C^\infty_0 ([0, T] \times
\Bbb{R}^3)$ with $\triangledown_x \cdot \psi = 0$ and for every
$0 \le t \le T$, we have
\begin{equation}
u(t)\psi(t) - u(0) \psi (0) - \int^t_0 u \cdot \partial_s \psi ds
= \int^t_0 \int_{\Bbb{R}^{3}} u \cdot \triangledown \psi \cdot u
dx ds 
\end{equation}
and $\triangledown_x \cdot u = 0$ in the sense of distributions.

We define the Littlewood-Paley energy flux $\Pi_j$ through a sphere
in frequency space of radius $2^j$ as follows.  For any divergence
free vector field $u \in L^2 (R^3)$ we define
\begin{equation}\label{}
S_ju = u * {\mathcal{F}}^{-1} (\psi (\cdot 2^{-j}))
\end{equation}
where $\psi(\xi)$ is a smooth nonnegative function supported in the
ball of radius one centered at the origin and such that $\psi (\xi)
= 1$ for $\xi \le 1/2$ and $\mathcal{F}$ is the Fourier transform. We then
define $\Pi_j$ as
\begin{equation}\label{}
\Pi_j = -\int_{\Bbb{R}^{3}} u \cdot \triangledown S^2_j u \cdot u \;
dx.
\end{equation}
Using the test function $S^2_j u$ in the weak formulation of the
Euler equations we obtain
\begin{equation}\label{}
\frac{1}{2} \frac{d}{dt} \|S_j u\|_2^2 = -\Pi_j 
\end{equation}
Cheskidov et al \cite{ccfs} prove that the Littlewood-Paley
energy flux of a divergence free vector field $u \in L^2$ satisfies
the following estimate:
\begin{equation}\label{}
|\Pi_j| \lesssim \sum\limits^\infty_{i = -1} 2^{-\frac{2}{3}|j-i|}
2^i \|u_i\|^3_3 
\end{equation}
where $u_j$ is the $j$-th Littlewood-Paley piece of $u$ defined by
$$
u_j = S_{j+1} u - S_j u.
$$
It follows from (1.19) that
\begin{equation}\label{}
\limsup_{j \to \infty} |\Pi_j| \lesssim C \limsup_{j \to
\infty} (2^j \|u_j \|^3_3). 
\end{equation}
Furthermore, an important feature of the bound (1.19) is that it is
quasi-local in the sense of rapid decay when $|j-i|$ is large.

We define the Besov space $B^{1/3}_{3, c_0}$ to be
the space of all tempered distributions $u \in \Bbb{R}^3$ for which
\begin{equation}\label{}
\limsup_{j \to \infty} 2^{j/3} \|u_j\|_3 = 0.
\end{equation}
Hence, from (1.18) and (1.20) we obtain the following result.

\emph{ Every weak solution $u$ to the Euler
equation on a time interval $[0, T]$ which satisfies
\begin{equation}\label{e:zerofluxcond}
\lim_{j \to \infty} \int^T_0 2^j \|u_j(t)\|^3_3 = 0 
\end{equation}
conserves energy on the entire interval $[0, T]$.  In particular,
energy is conserved for every solution in the class $L^3[0, T];
B^{1/3}_{3, c_0}) \cap C_w ([0, T]; L^2).$ }

In order to see more transparently the connection between (1.22) and
the smoothness 1/3 predicted by Onsager we rewrite (1.22) as
follows:
\begin{equation}\label{}
\lim_{|y|\to 0} \frac{1}{|y|} \int^T_0 \int_{\Bbb{R}^3} |u(x)
- u(x-y)|^3 dx dt = 0. 
\end{equation}
Hence the solution to the Euler equation needs to be a little
better than 1/3 Holder continuous in the space-time average to
ensure that energy is conserved.
We call the Besov space $B^{1/3}_{3, \infty}$ Onsager critical. This
is the space which contains distributions $u \in \Bbb{R}^3$ where
$\lim_{j \to \infty} 2^{j/3}\|u_j\|_3$ is finite, but not
necessarily zero.  This is a critical space in which energy
conservation for the Euler equation might be violated.

Applying the bound on the energy flux given by (1.19) to the
Navier-Stokes equations gives a sufficient condition for the energy
equality to hold, namely.

\emph{ Let $u^\nu \in C_w([0, T]; L^2
(\Bbb{R}^3)) \cap L^2 ([0, T]; H^1(\Bbb{R}^3))$ be a weak solution
to (1.6)-(1.7) with
\begin{equation}\label{}
\lim_{j \to \infty} \int^T_0 2^j \|u^\nu_j(t)\|^3_3 dt = 0.
\end{equation}
Then $u^\nu$ satisfies the energy equality
\begin{align}
\|u^\nu(t)\|^2_2 + 2\nu \int^t_0 \| \triangledown u^\nu(s) \|^2_2
ds
= \| u^\nu(0) \|^2_2 + 2 \int^t_0 f \cdot u^\nu(s) ds.
\end{align}
In particular, (1.25) holds if $u \in L^3([0, T]; H^{5/6})$.}

\subsection{Littlewood-Paley Framework for Intermittency}

Let $u^\nu$ be a Leray-Hopf weak solution to the Navier-Stokes
equations \eqref{e:nse} - \eqref{e:nsediv}.  We denote by $\langle \cdot \rangle$ the
long time average. We define the Littlewood-Paley energy spectrum of
$u^\nu$ by
\begin{equation}\label{e:LPspec}
E_{LP} (\k) = \frac{\langle\|u^\nu_j\|^2_2\rangle}{\k}
\end{equation}
for frequencies $\k \in [2^j, 2^{j+1}]$ and we define the mean energy
dissipation rate by
\begin{equation}\label{}
\varepsilon^\nu = \nu \langle \|\triangledown u^\nu\|^2_2\rangle.
\end{equation}
If a family of individual realizations $\{u^\nu\}_{\nu < \nu_0}$
verifies Kolmogorov's hypothesis that $\varepsilon^\nu \to
\varepsilon > 0$, then the locality of the flux which is exhibited
in the bound (1.19) suggests the following
\begin{equation}\label{e:limitbesov}
2^j \langle \| u^0_j \|^3_3\rangle \sim \varepsilon 
\end{equation}
for all $j$ sufficiently large. Here $u^0 = \lim_{\nu \to 0}
u^\nu$.  In other words, the limiting solution $u^0$ to the Euler
equation is ``on average" in the Onsager critical space $B^{1/3}_{3,
\infty}$.

Eyink \cite{eyink-besov} showed that $B^{1/3}_{3, \infty}$ is
consistent with the multi-fractal intermittency models of Frisch and
Parisi \cite{fp}.  Within the Littlewood-Paley framework we
can model the intermittency correction, (see \eqref{e:specinter} and \eqref{e:dissinter}) by
assuming the relationship between $\varepsilon$ and $\| u^0_j \|_3$
given in \eqref{e:limitbesov} and fixing the saturation level in Bernstein's
inequalities.  More precisely, in 3 dimensions we have
\begin{equation}\label{}
\|u_j\|_3 \lesssim 2^{j/2} \|u_j\|_2. 
\end{equation}
Assuming that the region of active turbulence is bounded, (say, on a
torus) we also have
\begin{equation}
\| u_j \|_2 \lesssim \|u_j\|_3. 
\end{equation}
Hence for $2^j \in [\k_f, \k_d]$ we write
\begin{equation}\label{e:saturlevel}
\|u_j\|_3 \sim 2^{cj} \|u_j\|_2 
\end{equation}
for some $c$ in the interval $[0, 1/2]$.  Then from \eqref{e:LPspec}, \eqref{e:limitbesov}
and \eqref{e:saturlevel} we recover the energy spectrum law
\begin{equation}
E_{LP} (\k) \sim \frac{\varepsilon^{2/3}}{\k^{5/3 + 2c}} 
\end{equation}
with $2^j$ being identified with $\k$.

The analogy between \eqref{e:specinter} and \eqref{e:saturlevel} requires that $D = 3-6c$. So,
the fully saturated Bernstein's inequality (i.e., $c = 1/2)$
corresponds to a uniform distribution of modes $u_j$ in each dyadic
shell and hence strong localization in space (i.e. $D=0$).  On the other
extreme, $c=0$ corresponds to a uniform distribution of $u$ in physical
space space (i.e. $D=3$) and localization in frequency space which
corresponds to the classical $K41$ model.

\section{A Continuous Model for the Energy Flux}

Although there is abundant empirical evidence for Kolmogorov's
hypothesis that $\lim_{\nu \to 0} \varepsilon^\nu \to \varepsilon
> 0$, this has not been rigorously proved for the Navier-Stokes to
Euler limit.  It is therefore of interest to examine simpler models
that retain some of the essential features of the fluid equations
and yet are tractable enough to allow a proof of Kolmogorov's
hypothesis. We now propose a PDE model for the turbulent energy
spectrum in frequency space.  We choose the scaling to include the
intermittency correction that we described in sections 1.2 and 1.3.
To motivate the model
we start with a fully local version of the flux given by the bound
in the inequality (1.20), namely
\begin{equation}\label{}
\Pi_j \sim 2^j \|u_j\|^3_3.
\end{equation}
We further assume that $c$ in \eqref{e:saturlevel} is independent of $j$. We thus obtain
\begin{equation} \label{e:flux1}
\Pi_j \sim 2^{(3c+1)j} \|u_j\|_2^3.
\end{equation}
We now make a further step from the discrete expression for flux \eqref{e:flux1} to a continuous one by looking at the energy density function $a(\k,t)$ defined by
$$
\| u(t)\|_2^2 = \int_0^\infty a^2(\k,t) d\k.
$$
Assuming that $a$ "does not vary much" in each dyadic shell, or disregarding energy density variations in each dyadic shell, we obtain
$$
\|u_j\|^2_2 \sim \int_{2^j}^{2^{j+1}} a^2(\k,t) d\k \sim \k a^2(\k,t).
$$
for all $\k \in [2^j,2^{j+1}]$. Thus, \eqref{e:flux1} becomes
\begin{equation} \label{e:flux2}
\Pi(\k) = \k^{3c + \frac{5}{2}} a^3(\k,t).
\end{equation}
Going back to \eqref{e:fluxnse} we assume that $\k_f = 1$, and the energy in the sub-inertial range in negligible.  So, with \eqref{e:flux2} at hand we can write the energy balance relation as follows
$$
\frac{1}{2} \left( \int_1^\k a^2(\ell,t) d\ell \right)_{t} = - \k^{3c + \frac{5}{2}} a^3(\k,t) - \nu \int_1^\k \ell^2 a^2(\ell,t) d\ell.
$$
Differentiating in $\k$ and cancelling $a$ on both sides we obtain the following PDE
\begin{equation}\label{e:model1}
a_{t} = -(\frac{5}{2} + 3c ) \k^{\frac{3}{2} + 3c} a^{2} - 3\k^{\frac{5}{2} + 3c} aa_{\k} - \nu \k^{2} a
\end{equation}
We supplement this equation  with the boundary condition 
\begin{equation}\label{e:bc}
a(1,t)=\e^{1/3}.
\end{equation}
Here $\e$ represents the energy input rate coming from an external force. We thus disregard any particular detail of energy production and simply model it with our boundary condition. The input rate $\e$ will subsequently be shown equal to the energy dissipation rate, hence the notation.
 Equation \eqref{e:model1} can be easily simplified by rescaling the energy density $a$ to $b = \k^{\frac{5}{6} + c} a$. Thus, \eqref{e:model1} becomes the following equation
\begin{equation}\label{e:model2}
b_{t} = -3 \k^{\a} bb_{k} - \nu \k^{2} b, \quad \k \geq 0,
\end{equation}
where $\a = \frac{5}{3} + 2c$. The appropriate range of $\a$ is $[\frac{5}{3}, \frac{8}{3}]$, which exactly corresponds to the classical range of the energy density power laws with the spatial intermittency correction.

\begin{remark}

Equation \eqref{e:model2} is a variant of the much studied Burgers equation.
The one dimensional Burgers equation can be viewed as the most basic
nonlinear PDE that has the bilinear structure of the nonlinearity of
the Euler and Navier-Stokes equations. It can be invoked as a model
for one dimensional compressible fluids. However there is no clear
physical basis for using Burgers equation in physical space as a
model for turbulence. On the other hand, as we have argued, the
locality of the energy flux manifested using Littlewood-Paley
theory, motivates \eqref{e:model2} as a PDE model for the turbulent cascade in
frequency space.
\end{remark}

\begin{remark}
In the past few decades a number of ``toy models" for turbulence
have been studied to test Kolmogorov's theory. 
In particular, the derivation of the classical Desnyanskiy-Novikov discrete model, \cite{dn}, follows a similar path. The flux there is modeled by taking $\Pi_{j} = 2^{dj}a_{j}a_{j+1}$, where $a_{j}^{2}$ represents the total energy in the $j$-th dyadic shell, while $d$ is an intermittency parameter with the appropriate range of values. The model is thus an infinite system of ODEs given by
\begin{equation} \label{e:dyadic}
\frac{d}{dt} a_j + \nu 2^{2j}a_j - 2^{d(j-1)}a_{j-1}^2 + 2^{dj}
a_j a_{j+1}=0, \qquad j=0,1,2,\dots,
\end{equation}
where $a_{-1}=0$.
This model, as well as its inviscid versions, has been extensively
studied by Katz and Pavlovic \cite{kp}, Cheskidov and Friedlander \cite{cf}, Kiselev and Zlatos \cite{kz}, and others.
To our knowledge the PDE model
we present is the first continuous model. 
\end{remark}

In Section \ref{s:inviscid} we examine the inviscid $(\nu = 0)$ form of \eqref{e:model1} with
boundary condition \eqref{e:bc}. We prove
that there is a unique fixed point which is a global attractor. Moreover every
solution reaches it in finite time. The inviscid equation
exhibits anomalous dissipation and the average energy spectrum has
the power law $\e^{2/3} \kappa^{-\alpha}$.

In Section \ref{s:viscous} we turn to the viscous model $\nu > 0$ \eqref{e:model1}, which is in essence a damped Burgers equation. Again the PDE
has a unique fixed point and this converges to the inviscid fixed
point as $\nu \to 0$. The viscous fixed point reproduces
Kolmogorov's energy density spectrum in the inertial range and it becomes zero identically  after the dissipation wave number $\kappa_d$. We further consider the Leray regularization of equation \eqref{e:model1}. We show that all bounded solutions of the regularized equation converge pointwise and in the metric of $L^2$-space to a fixed point which in turn converges to the fixed point of equation \eqref{e:model1} in the limit as the regularization parameter goes to zero. 
The average dissipation rate for the viscous system converges
to the anomalous dissipation rate for the inviscid system giving an example of the dissipation anomaly predicted by Kolmogorov.

\section{Inviscid case}\label{s:inviscid}
In this section we study the inviscid version of the model \eqref{e:model1}
\begin{equation}\label{e:inviscid}
\left\{
\begin{split}
&a_{t} = -3a \k^{\alpha} \frac{\partial}{\partial \k}(\k^{\frac{\a}{2}}a), \qquad \k>1,\\
&a(1,t)=\e^{1/3},\\
&a(\k,0)=a_0(\k) \geq0,
\end{split}
\right.
\end{equation}
where $\alpha \in [5/3,8/3]$.
Note that the energy equality
\begin{equation}
\frac{d}{dt}\frac{1}{2}\int_1^{\infty} a(\k,t)^2 \, d\k = \e
\end{equation}
is satisfied on some interval $t\in(0,T)$ provided the solution satisfies
the following smoothness condition (cf. \eqref{e:zerofluxcond}):
\begin{equation} \label{encon}
\lim_{\k \to \infty} \int^T_0 \k^{3\alpha/2} a^3 \, dt = 0,
\end{equation}
which is an analog to \eqref{e:zerofluxcond}. So, $\e>0$ represents the energy input rate in this model.
The unique fixed point of \eqref{e:inviscid} is given by 
\begin{equation}
A^0(\k) = \e^{1/3}\k^{-\alpha/2}
\end{equation}
We note that the fixed point does not satisfy
\eqref{encon}. Moreover, it does not satisfy the energy equality since
\begin{equation}
\frac{d}{dt}\frac{1}{2}\int_1^{\infty} A^0(\k)^2 \, d\k = 0 \ne \e.
\end{equation}
The anomalous
energy dissipation rate is the difference between the energy input rate
and the time derivative of the total energy. Thus using (3.2) and (3.5)
we observe that the anomalous dissipation rate for the fixed point 
is exactly the energy
input rate $\e$. We will show that this also holds for every other solution
asymptotically in time. In order to do this we will prove that $A^0$ is a global attractor.

We use the following change of variables:
\[
\xi=\k^{-1/\g}, \qquad
w(\xi,t)=\e^{-1/3}\xi^{-\frac{\g\a}{2}}a(\xi^{-\g},{\textstyle \frac{1}{3}}\e^{-1/3}\g t),
\]
where $\g=\frac{1}{\a-1}$. Then \eqref{e:inviscid} reduces to
Burgers equation
\begin{equation}
w_t = w w_{\xi}, \qquad 0<\xi<1
\end{equation}
with $w(1,t)=1$, $w(\xi, 0)=w_0(\xi) \geq0$. We extend it to Burgers equation on the whole real line
\begin{equation} \label{e:burgers-extended}
\left\{
\begin{split}
&w_t = w w_{\xi}, \qquad \xi \in \mathbb{R},\\
&w(\xi,0) = \tilde{w}_0(\xi),
\end{split}
\right.
\end{equation}
where
\[
\tilde{w}_0(\xi)=\left\{
\begin{split}
&0, \qquad &\xi<0,\\
&w_0(\xi), \qquad &0\leq \xi \leq 1,\\
&1, \qquad &\xi> 1.
\end{split}
\right.
\]

The weak solution to \eqref{e:burgers-extended} is expressed using
the Lax-Oleinik formula. Let
\[
h(y)=\int_0^y \tilde{w}_0(\xi) \, d\xi.
\]
For all $t$ and for all but at most countably many $\xi \in \mathbb{R}$,
there exists a unique  $y_*(\xi,t)$, such that
\[
\min_y\left\{f(y) \right\} = f(y_*).
\]
where
\[
f(y)=\frac{(\xi-y)^2}{2t}-h(y).
\]
Then
\begin{equation} \label{e:L-A}
w(\xi,t) = \frac{y_*(\xi,t)-\xi}{t}
\end{equation}
is the weak solution to \eqref{e:burgers-extended}.
Given $\xi \in [0,1]$, let $t>2(1-\xi)$. We will show that
$y_*(\xi,t)=\xi+t$ and consequently $w(\xi,t)=1$. Indeed, let $y=\xi+s$.
First, consider the interval $s\geq 1-\xi$. Then
\[
f(y) = \frac{s^2}{2t} - \int_0^1 w_0(\xi) \, d\xi -\xi-s+1 > f(\xi+t) 
\]
provided $s\ne t$. Since $t>2(1-\xi)$, it follows that
\[
f(\xi+t) < - \int_0^1 w_0(\xi) \, d\xi,
\]
and hence
on the interval $s<1-\xi$ we also have
\[
f(y) \geq \frac{s^2}{2t} - \int_0^1 w_0(\xi)  \, d\xi \geq
 - \int_0^1 w_0(\xi)  \, d\xi> f(\xi+t).
\]
Therefore, \eqref{e:L-A} implies that $w(\xi,t) =1$ for $\xi \in (0,1)$, $t \geq 2$ and hence, returning to the original variables,
\begin{equation}
a(\k, t) = A^0(\k), \qquad  t\geq {\textstyle \frac{2}{3}}\e^{-1/3}\g.
\end{equation}
Hence the average energy spectrum for solutions to the inviscid model
is about $\e^{2/3} \k^{-\alpha}$.

\section{Viscous case}\label{s:viscous}

In this section we study the viscous model \eqref{e:model1}
\begin{equation}\label{e:viscous}
\left\{
\begin{split}
&a_{t} = -3a \k^{\alpha} \frac{\partial}{\partial \k}(\k^{\frac{\a}{2}}a)
-\nu \k^2a,
 \qquad \k>1,\\
&a(1,t)=\e^{1/3},\\
&a(\k,0)=a_0(\k), a_0 \in L^2, a_0\geq0,
\end{split}
\right.
\end{equation}
where $\alpha \in [5/3,8/3]$.
There exists a unique fixed point to \eqref{e:viscous} given by 
\begin{equation} \label{fixed-point-visA}
A^\nu(\k)=\left\{
\begin{aligned}
& \k^{-\a/2} \left[ \e^{1/3}+ \frac{\nu}{3(3-\a)}(1- \k ^{3-\a}) \right],  &1 \leq \k \leq \k_d,\\
&0,  & \k > \k_d,
\end{aligned}
\right.
\end{equation}
where $\k_d$ is Kolmogorov's dissipation wavenumber described in
Section 1.1. For the model it is explicit and given by
\begin{equation}
\k_d= \left[   1+ \frac{3(3-\a) \e^{1/3}}{\nu}  \right]^{\frac{1}{3-\a}}.
\end{equation}
To see the parallel with the classical expressions for $\k_d$ we note that for $\nu$ small one has
\[
\k_d \sim \left(\frac{\e}{\nu^3}\right)^{\frac{1}{4}}, \qquad \text{for} \qquad \alpha =\frac{5}{3},
\]
and
\[
\k_d \sim \frac{\e}{\nu^3}, \qquad \text{for} \qquad \alpha =\frac{8}{3}.
\]

In the limit of vanishing viscosity we immediately obtain from (4.2) and (3.4)
that
\begin{equation}
A^\nu(\k) \ra A^0(\k),
\end{equation}
uniformly on any finite interval $[1,\k_0]$ as $\nu\to 0$. With a little more effort we can show that the convergence also takes place in $L^2([1,\infty))$. Indeed,
\begin{multline*}
\int_1^\infty |A^\nu(\k) - A^0(\k)|^2 d\k \leq \e^{2/3} \int_{\k_d}^\infty \k^{-\a} d\k + \\ \frac{\nu^2}{9(3-\a)^2} \int_1^{\k_d} k^{-\a}( 1-\k^{3-\a}  )^2 d\k.
\end{multline*}
Since $\k_d \ra \infty$ we see that the first integral vanishes as $\nu \ra 0$. The second integral behaves like $\nu^2$ for $\a <7/3$, like $\nu^2 \log(\nu)$ for $\a = 7/3$, and like $\nu^{(\a-1)/(3-\a)}$ for $7/3< \a \leq 8/3$. So, within our range of $\a$ the second integral vanishes too as $\nu\ra0$.

In order to study the time dependent solutions to the viscous system \eqref{e:viscous} we utilize the same change of variables as in the previous section
\[
\xi=\k^{-1/\g}, \qquad
w(\xi,t)=\e^{-1/3}\xi^{-\frac{\g\a}{2}}a(\xi^{-\g},{\textstyle \frac{1}{3}} \e^{-1/3}\g t),
\]
where $\g=\frac{1}{\a-1}$. Then \eqref{e:viscous} reduces to the
damped Burgers equation
\begin{equation} \label{eq:dumpedBurgers}
w_t = w w_{\xi} - \mu \xi^{-2\gamma} w, \qquad 0<\xi<1
\end{equation}
with $w(1,t)=1$, $w(\xi, 0)=w_0(\xi) \geq0$, $\int_0^1 \xi^2 w_0^2(\xi) \, d\xi <\infty$. Here $\mu= \frac{1}{3}\nu \e^{-1/3}\gamma$.
The equation has a unique fixed point
\begin{equation} \label{fixed-point-vis}
W(\xi)=\left\{
\begin{aligned}
&1+\frac{\mu}{2\gamma-1}(1-\xi^{1-2\gamma}),  &\xi_d< \xi \leq 1,\\
&0,  &0 \leq \xi \leq \xi_d,
\end{aligned}
\right.
\end{equation}
where
\begin{equation}
\xi_d = \left[1+\frac{2\gamma-1}{\mu}\right]^{\frac{1}{1-2\gamma}}.
\end{equation}
Note that $\xi_d \to 0$ as $\mu \to 0$. Note that $\k_d =\xi_d^{-\gamma}$.

Since there are discontinuous solutions to the damped Burgers equation
\eqref{eq:dumpedBurgers}, we now 
consider a Leray-type regularization of the equation. Such regularizations
have been used to approximate weak solutions to the Burgers equations
\cite{BF, NM}. Consider the following regularized equation:
\begin{equation} \label{eq:reg}
w_t = v_\delta w_\xi - \mu \xi^{-2\gamma} w, \qquad v_\delta = w*\phi_\delta, \qquad 0<\xi<1, t>0,
\end{equation}
with the boundary conditions $w(1,t)=1$, $w(0,t)=0$.
Here $\phi_\delta(\xi)=\frac{1}{\delta}\phi(\xi/\delta)$, where $\phi(\xi)$ is smooth,
nonnegative, and such that $\int \phi(\xi) \, d\xi =1$, $\supp \phi = (-1,0)$.
First, note that there exists a unique fixed point $w=W_\delta(\xi)$ to
\eqref{eq:reg}, which is smooth, nonnegative, monotonically increasing,
and with $W_\delta(0)=0$.
Now consider characteristics $\eta_\d(t)$, which are solutions to
\begin{equation} \label{e:char}
\ddt \eta_\d(t)=-v_\d(\eta_\d(t),t).
\end{equation}
It is easy to see that characteristics do not intersect on $(0,1)$. Indeed,
take $\epsilon \in (0,1)$ and consider
two solutions $\eta'_\d(t)$, $\eta''_\d(t)$ to \eqref{e:char}. As long as
they belong to $(\epsilon, 1)$ we have
\[
\begin{split}
\ddt|\eta'_\d(t)-\eta''_\d(t)| &\lesssim
\delta^{-3/2}|\eta'_\d(t)-\eta''_\d(t)|
\left(\int_\epsilon^1  w_\d^2(\xi,t)\, d\xi \right)^{1/2}\\
&\lesssim
\delta^{-3/2} |\eta'_\d(t)-\eta''_\d(t)|
\frac{1}{\epsilon}\left(\int_\epsilon^1 \xi^2 w_\d^2(\xi,t)\, d\xi \right)^{1/2}.
\end{split}
\]
Since $\int_0^1 \xi^2 w_\delta^2(\xi,t)\, d\xi$ is non-increasing, the characteristics
do not intersect.

Along characteristics we have
\[
\ddt w_\d(\eta_\d(t),t) = -\mu\eta_\d^{-2\gamma} w_\d(\eta_\d(t),t).
\]
Moreover, along every characteristic curve that starts from the boundary
$\xi=1$ we have that $w_\d$ is equal to the fixed point $W_\d$, i.e.,
\[
w_\d(\eta_\d(t),t) = W_\delta(\eta_\d(t)),
\]
provided $\eta_\d(t_0)=1$ for some $t_0\geq 0$. Now consider the characteristic
curve $\eta_\d^0(t)$ with $\eta_\d^0(0)=1$. Note that $\eta_\d^0(t)$ is decreasing,
positive, $\eta_\d^0(t) \to 0$ as $t \to \infty$, and
\begin{equation}
w_\d(\xi,t) = W_\d(\xi),
\end{equation}
provided $\xi \geq \eta_\d^0(\xi,t)$.
Note also that $\eta^0_{\d_1}(t) \leq \eta^0_{\d_2}(t)$ provided $\d_1\leq \d_2$. Then consider
the particular characteristic curve $\eta^0_{1}(t)$ corresponding to $\d =1$. Let $T(\xi)$ be such that
$\eta^0_1(T(\xi)) =\xi$. Let $\xi_0 \in (0,1)$. Then for all $\d \leq 1$, $\xi \in[\xi_0, 1]$, $t\geq T(\xi_0)$ we have
\begin{equation}
w_\d(\xi,t) = W_\d(\xi).
\end{equation}
So the fixed point of the regularize equation $W_\delta$ is a
global attractor.
Furthermore, as $\delta \to 0$ the fixed point $W_\delta$ converges
to the fixed point $W$ given by \eqref{fixed-point-vis}.

We now return to equation \eqref{e:viscous} for $a(\xi, t)$.
We will study the energy dissipation in the limit of
vanishing viscosity. For convenience, solutions to the
model \eqref{e:viscous} with viscosity $\nu \geq 0$ will be denoted by
$a^\nu(t)$. The fixed point (which is unique
in both viscous and inviscid cases) will be denoted by $A^\nu$.
We assume that $a^\nu(\cdot,0) \in L^2$ and
$a^\nu(t) \to A^\nu$ in $L^2$. We denote the energy and enstrophy
norms of $a^\nu$ by
\[
|a^\nu(t)|:=\int_1^\infty a^\nu(\k,t)^2 \, d\k, \qquad \|a^\nu(t)\|:=\int_1^\infty \k^2a^\nu(\k,t)^2 \, d\k.
\] 

Due to the energy inequality we have
\[
\frac{1}{2T}|a^\nu(T)|^2 - \frac{1}{2T}|a^\nu(0)|^2 \leq -\nu \frac{1}{T}
\int_0^{T} \|a^\nu(t)\|^2 \, dt + \e.
\]
Hence,
\[
\limsup_{T \to \infty} \frac{1}{T} \int_0^T \nu \|a^\nu(t)\|^2 \, dt \leq  \e.
\]
On the other hand, note that the fixed point $A^\nu$
satisfies the energy equality
\[
\nu \|A^\nu\|^2 = \e.
\]
Now for any $\eta>0$, there exists $N$, such that
\[
\nu \int_{1}^N \k^2 A^\nu(\k)^2 \, d\k \geq \nu \|A^\nu\|^2 - \eta.
\]
Since $a^\nu(t) \to A^\nu$ in $L^2$,  we have
\[
\liminf_{T \to \infty} \frac{1}{T} \int_0^T  \int_1^N \nu \k^2 a^\nu(\k,t)^2\, d\k \, dt \geq \nu \int_1^N \k^2 A^\nu(\k)^2 \, d\k\geq \nu \|A^\nu\|^2 - \eta.
\]
Therefore,
\[
\begin{split}
\liminf_{T \to \infty} \frac{1}{T} \int_0^T \nu \|a^\nu(t)\|^2 \, dt \geq
\nu\|A^\nu\|^2 = \e.
\end{split}
\]
Then we obtain
\[
\e:=\lim_{\nu\to 0} \lim_{T \to \infty} \frac{1}{T} \int_0^T \nu \|a^\nu(t)\|^2 \, dt >0.
\]
i.e., in the limit of vanishing viscosity the energy dissipation rate is positive and equal to the energy input rate or anomalous dissipation rate for the 
inviscid model.


\begin{thebibliography}{10}

\bibitem{BF}. S. Bhat and R. C. Fetecau. A Hamiltonian regularization of the
Burgers equation, {\it J. Nonlinear Sci.}, 16(6):615Ð638, 2006.

\bibitem{ccfs}
A.~Cheskidov, P.~Constantin, S.~Friedlander, and R.~Shvydkoy.
\newblock Energy conservation and {O}nsager's conjecture for the {E}uler
  equations.
\newblock {\em Nonlinearity}, 21(6):1233--1252, 2008.

\bibitem{cf}
Alexey Cheskidov and Susan Friedlander.
\newblock The vanishing viscosity limit for a dyadic model.
\newblock {\em Phys. D}, 238(8):783--787, 2009.

\bibitem{cet}
Peter Constantin, Weinan E, and Edriss~S. Titi.
\newblock Onsager's conjecture on the energy conservation for solutions of
  {E}uler's equation.
\newblock {\em Comm. Math. Phys.}, 165(1):207--209, 1994.

\bibitem{ds}
Camillo De~Lellis and L{\'a}szl{\'o} Sz{\'e}kelyhidi, Jr.
\newblock The {E}uler equations as a differential inclusion.
\newblock {\em Ann. of Math. (2)}, 170(3):1417--1436, 2009.

\bibitem{dn}
V.N. Desnyansky and E.A. Novikov.
\newblock The evolution of turbulence spectra to the similarity regime.
\newblock {\em Izv. Akad. Nauk SSSR Fiz. Atmos. Okeana}, 10:127--136, 1974.

\bibitem{dr}
Jean Duchon and Raoul Robert.
\newblock Inertial energy dissipation for weak solutions of incompressible
  {E}uler and {N}avier-{S}tokes equations.
\newblock {\em Nonlinearity}, 13(1):249--255, 2000.

\bibitem{e}
Gregory~L. Eyink.
\newblock Energy dissipation without viscosity in ideal hydrodynamics. {I}.
  {F}ourier analysis and local energy transfer.
\newblock {\em Phys. D}, 78(3-4):222--240, 1994.

\bibitem{eyink-besov}
Gregory~L. Eyink.
\newblock Besov spaces and the multifractal hypothesis.
\newblock {\em J. Statist. Phys.}, 78(1-2):353--375, 1995.
\newblock Papers dedicated to the memory of Lars Onsager.

\bibitem{es}
Gregory~L. Eyink and Katepalli~R. Sreenivasan.
\newblock Onsager and the theory of hydrodynamic turbulence.
\newblock {\em Rev. Modern Phys.}, 78(1):87--135, 2006.

\bibitem{fp}
U.~Frisch and G.~Parisi.
\newblock On the singularity structure of fully developed turbulence.
\newblock In M.~Ghil, R.~Benzi, and G.~Parisi, editors, {\em Turbulence and
  predictability in geophysical fluid dynamics and climate dynamics}, Proc.
  International Summer School of Physics "Enrico Fermi", pages 84--87.
  Amsterdam: North-Holland, 1985.

\bibitem{frisch}
Uriel Frisch.
\newblock {\em Turbulence}.
\newblock Cambridge University Press, Cambridge, 1995.
\newblock The legacy of A. N. Kolmogorov.

\bibitem{kp}
Nets~Hawk Katz and Nata{\v{s}}a Pavlovi{\'c}.
\newblock Finite time blow-up for a dyadic model of the {E}uler equations.
\newblock {\em Trans. Amer. Math. Soc.}, 357(2):695--708 (electronic), 2005.

\bibitem{kz}
A. Kiselev and A. Zlato{\v{s}}.
\newblock On discrete models of the {E}uler equation.
\newblock {\em Int. Math. Res. Not.}, (38):2315--2339, 2005.

\bibitem{K41}
A.~N. Kolmogorov.
\newblock The local structure of turbulence in incompressible viscous fluids at
  very large reynolds numbers.
\newblock {\em Dokl. Akad. Nauk. SSSR}, 30:301--305, 1941.


\bibitem{NM} G. Norgard and K. Mohseni.
On the Convergence of the Convectively Filtered 
Burgers Equation to the Entropy Solution of the 
Inviscid Burgers Equation, arXiv:0805.2176v5.

\bibitem{onsager}
L.~Onsager.
\newblock Statistical hydrodynamics.
\newblock {\em Nuovo Cimento (9)}, 6(Supplemento, 2(Convegno Internazionale di
  Meccanica Statistica)):279--287, 1949.

\bibitem{scheffer}
V. Scheffer.
\newblock An inviscid flow with compact support in space-time.
\newblock {\em J. Geom. Anal.}, 3(4):343--401, 1993.

\bibitem{shn}
A.~Shnirelman.
\newblock On the nonuniqueness of weak solution of the {E}uler equation.
\newblock {\em Comm. Pure Appl. Math.}, 50(12):1261--1286, 1997.

\end{thebibliography}

\end{document}